\documentclass{amsart}

\usepackage{amsmath, amsfonts, amssymb, amscd}

\newtheorem{theo}{Theorem}[section]
\newtheorem{lem}[theo]{Lemma}

\newtheorem{rem}[theo]{Remark}

\newtheorem{definition}[theo]{Definition}
\newenvironment{pf}{\noindent{\it Proof. }}{$\square$\par\medskip}

\newcommand{\R}{{\mathbb R}}

\newcommand{\g}{{\mathfrak g}}

\newcommand{\J}{{\mathbb J}}

\newcommand{\D}{{\mathcal D}}
\newcommand{\N}{{\mathcal N}}
\renewcommand{\L}{{\mathcal L}}
\newcommand{\BL}{{\mathbb L}}

\newcommand{\U}{{\mathcal U}}
\newcommand{\V}{{\mathcal V}}

\renewcommand{\=}{\overset{\text{def}}{=}}

\def\sideremark#1{\ifvmode\leavevmode\fi\vadjust{% The remark
\vbox to0pt{\hbox to 0pt{\hskip\hsize\hskip1em% will appear only
\vbox{\hsize3cm\tiny\raggedright\pretolerance10000% on the side
\noindent #1\hfill}\hss}\vbox to8pt{\vfil}\vss}}}% in 3cm

\title[Total reality of conormal bundles]
{Total reality \\
of conormal bundles of  hypersurfaces \\in almost complex manifolds}

\author{A.Spiro}

\date{\today}


\begin{document}

\begin{abstract}
A generalization to the almost complex setting of a well-known
result by S. Webster is given. Namely,  we prove that  if $\Gamma$ is a strongly 
pseudoconvex hypersurface in an almost complex manifold $(M, J)$, then 
the  conormal bundle of $\Gamma$ is a totally real submanifold of $(T^*M, \J)$, where $\J$
is the lifted almost complex structure  on $T^*M$ defined by  Ishihara and Yano.

\end{abstract}

\maketitle

\null \vspace*{-.25in}

\section{Introduction}
\bigskip
Let $(M, J)$ be an almost complex manifold of real dimension $2n$ and $\J$ the associated
almost complex structure on $T^*M$ defined by Ishihara and Yano in \cite{IY} (see definition in \S 2, below). Consider  a smooth real hypersurface $\Gamma \subset M$ and  denote by 
$\N(\Gamma) \subset T^*M$ the   conormal bundle of $\Gamma$, i.e. the
 submanifold of $T^*M$ defined by 
$$\N(\Gamma) = \bigcup_{x\in \Gamma}  \N(\Gamma)_x\ ,\qquad \text{where}\ \  \N(\Gamma)_x  \= \{\ \alpha \in T^*_xM \ :\ \alpha|_{T_x\Gamma} = 0\ \}\ .$$
The almost complex
structure $J$ of $M$ induces on $\Gamma$ a possibly non-integrable CR structure
$(\D, J)$, i.e. a distribution $\D$ given by the $J$-invariant subspaces of the tangent 
spaces of $\Gamma$ and the complex structures $J_x \= J|_{\D_x}$ on the subspaces of the 
distribution $\D$. \par
In this short note we prove  that  {\it if  the (possibly non-integrable) CR structure  $( \D, J)$ on $\Gamma$ 
 is strongly pseudoconvex 
or  if the almost complex structure $J$ is integrable and 
$\D$ is a contact distribution, then the complement of the zero section 
$\N(\Gamma)\setminus \{\text{zero section}\}$ 
is a totally real 
submanifold of $(T^*M, \J)$\/}Ê(see Theorem \ref{maintheorem}, below).\par
\medskip
This  fact was first proved by S. Webster in \cite{We}  for strongly  pseudoconvex hypersurfaces in a complex manifold $M$.  Later,  always assuming that $M$ is  a complex manifold, 
 the result  was generalized for the conormal bundles of 
 Levi non-degenerate hypersurfaces and for Levi non-degenerate 
 submanifolds  of codimension higher then one  by A. Tumanov in (\cite{Tu}). 
 Another proof for the Levi non-degenerate hypersurfaces in complex manifold was given
  by Z. M. Balogh and C. Leuenberger in (\cite{BL}). At the best of our knowledge, the result under the weaker 
assumption that  $J$ is a possibly non-integrable almost complex structure  was not previously
known. 
 \par
We have to point out  that requiring
the   conditions $J$ integrable  and  $\D$   contact is equivalent to assume the
 Levi non-degeneracy  of  the hypersurface $\Gamma$  in the complex manifold $(M,J)$. This means that 
 the second part of our claim is just equivalent to   Tumanov's generalization for hypersurfaces
 in complex manifolds. \par
\medskip
The question whether Webster's result could be valid also in the almost complex setting
was asked to the author by A. Sukhov, being an interesting  problem related to the theory 
of $J$-holomorphic discs attached to the boundary of strongly pseudoconvex domains. 
Indeed, as for the classical case (see e.g. \cite{We}), 
our result should turn out to be quite useful for proving  boundary  
regularity properties  of $J$-biholomorphisms  between bounded domains in almost complex manifolds.
Our result has been also used by H. Gaussier and A. Sukhov in their recent 
paper \cite{GS}.
\par
The author is grateful to A. Sukhov for telling him the problem and for useful discussions on the subject.\par
\medskip
We tried to make the paper as much as possible complete and self-contained:  After giving all needed  preliminaries
in \S 2  (e.g.  presentation of   Ishihara and Yano's lifted 
almost complex structure in  ``coordinate-free notation", definition of pseudoconvexity and Levi-forms for  hypersurfaces  in almost complex manifolds, etc.) the proof of our main result is given in \S 3.  \par
\bigskip
 \bigskip
\section{Basic definitions and preliminaries}
\bigskip
\subsection{Lift of an almost complex structure to the cotangent bundle}
Let $M$ be a $2n$-dimensional manifold endowed with an almost 
complex structure $J$ and  let $\pi: T^*M \to M$ be its cotangent bundle.
As we mentioned in the Introduction,  Ishihara and Yano  proved in   \cite{IY}  that 
{\it there exists a natural almost complex structure $\mathbb J$ on $T^*M$, such that: 
\begin{itemize}
\item[a)]  it is a "lift" of $J$ on $T^*M$, i.e. 
$\pi_*(\J V) = J \pi_*(V)$ for any vector field $V\in T(T^*M)$; 
\item[b)] it is 
invariant w.r.t. the lifted action on $T^*M$ of any $J$-preserving diffeomorphism $f: M \to M$.
\end{itemize}
}
In order to define such lifted almost complex structure on $T^*M$, we first have 
to introduce a few objects.\par
\medskip
First of all, let us  denote by $\theta$ the so-called {\it tautological 1-form of 
$T^*M$}, 
i. e. the 1-form on $T^*M$ defined by
$$\theta_\alpha(V) \= \alpha(\pi_*(V))$$
for any $\alpha\in T^*M$, $V\in T_\alpha(T^*M)$. It is known that $\omega = d\theta$ is a symplectic form on $T^*M$, which is called
{\it canonical symplectic form of $T^*M$\/}. \par
\medskip
We may also consider
the inverse tensor field $\omega^{-1} \in \Lambda^2 T(T^* M)$, i.e. the tensor 
field of type $(2,0)$ so that, 
for any $\alpha\in T^* M$, $V\in T_\alpha(T^* M)$
$$\omega^{-1}|_{\alpha}(\omega_\alpha(V, \cdot), \cdot) = V\ .\eqno(2.1)$$
It is clear that, in 
any given basis, the  components of $\omega^{-1}$ 
are the entries of  the inverse 
of the matrix  given by the components 
of $\omega$. It is also immediate to check  that for any 
$\alpha\in T^* M$ and  $A \in T^*_\alpha(T^* M)$
$$\omega_\alpha(\omega^{-1}|_{\alpha}(A, \cdot), \cdot) = A\ . \eqno(2.2)$$
\par
\medskip
Let us now consider 
the Nijenhuis tensor $N^J$ of $J$, i.e. the tensor field of type $(1,2)$ on $M$, 
 defined   as 
follows: for any given pair of vectors  $v,w \in T_xM$, let  us denote by 
$X^{(v)}$ and $X^{(w)}$ two 
vector fields such that $X^{(v)}|_x = v$ and $X^{(w)}|_x = w$; then $N^J_x(v,w)$ is defined by 
$$N^J(v,w) \= [JX^{(v)}, JX^{(w)}]_x - [X^{(v)},X^{(w)}]_x - J([X^{(v)},JX^{(w)}]_x + [JX^{(v)}, X^{(w)}]_x)\ .\eqno(2.3)$$
It is easily seen that   $N^J_x(v,w)$ is independent of the choice of the vector fields $X^{(v)}$
and $X^{(w)}$. The relevance of the tensor field $N^J$ is determined by the celebrated theorem of Newlander and Niremberg (\cite{NN}), for which
 {\it an almost complex structure $J$ is an (integrable) complex structure 
  if and only if $N^J \equiv 0$\/}.\par
\medskip
Now, with the help of  $N^J$ and of the symplectic form $\omega$, 
 we  may define the following tensor field $g^J$  of type $(0,2)$ 
 on $T^*M$:  for any $\alpha \in T^*M$ and  
$V,W \in T_\alpha(T^*M)$,  we set 
$$g^J_\alpha(V,W) \= \frac 1 2 \alpha(N^J(\pi_*(V),J\pi_*(W)) = $$
$$ = \frac 1 2  \alpha\left( - [JX^{(v)}, X^{(w)}]_x - [X^{(v)},J X^{(w)}]_x  + J[X^{(v)}, X^{(w)}]_x - J  [JX^{(v)}, J X^{(w)}]_x\right)\ ,\eqno(2.4)$$
where, as before,   $X^{(v)}$ and $X^{(w)}$ are two 
vector fields on $M$ such that $X^{(v)}|_x = v =  \pi_*(V) $ and $X^{(w)}|_x = w =  \pi_*(W)$. 
From  
definitions, it is clear that $g^J_\alpha(V,W)$ is skew-symmetric w.r.t $V$ and $W$ and 
that, if $\pi_*(W) = J \pi_*(V)$ then $g^J_\alpha(V,W) = 0$.\par
\medskip
The last necessary ingredient to define Ishihara and Yano's
  almost complex structure $\J$ on $T^*M$  is 
 the fiber preserving 
diffeomorphism $\hat J: T^*M \to T^*M$,  given by the  map which associates to any 1-form $\alpha \in T^*_x M$, $x\in M$, 
the 1-form 
$$\hat J(\alpha) \= \alpha \circ J\in T^*_x M\ .$$ 
\par
\medskip
Now we can give  the definition of lifted almost complex structure $\J$.\par
\begin{definition}{\rm Let   $\varpi^J$ be the tensor field of type 
$(0,2)$ on $T^*M$ defined by 
$$\varpi^J \= \hat J^*\omega + g^J = d(\hat J^* \theta) + g^J\ .\eqno(2.5)$$
We call  {\it  lifted almost complex structure on  $T^*M$ associated with $J$\/} 
the tensor field $\J$ of type $(1,1)$  defined by 
$$\J(v) \= \omega^{-1}(\varpi^J(V, \cdot), \cdot)\ , \qquad \text{for any}\ V\in T_\alpha(T^*M)\ .  \eqno(2.6) $$}
\end{definition}
\bigskip
It is proved in \cite{IY} that (2.6) does  define an almost complex structure on $T^*M$. Moreover, 
being $\hat J^*\omega$ and $g^N$ invariant under any lifted action of a  biholomorphism of $(M,J)$,  
it is immediate to realize  that the requirement b) of lifted almost complex structures holds.  
To check that  also  the condition a) is satisfied,  it suffices 
to write down  the explicit expressions of the components $\J$ in some 
coordinate basis.\par
\medskip
For this,  consider a coordinate chart 
$$\xi = (x^1, \dots, x^{2n}): \U \subset M \to \R^{2n}$$
and an associated chart 
$$\hat \xi = (x^1, \dots x^{2n}, p_1, \dots, p_{2n}): \pi^{-1}(\U) \subset T^*M \to \R^{4n}\ ,$$
where we denote by $p_i$'s the components 
 of the forms $\alpha = p_i dx^i \in \pi^{-1}(\U)$ w.r.t. the coordinate basis $dx^1$,  \dots, $dx^n$.
 Now, if we fix  $x \in \U \subset M$ and 
$\alpha = p_a dx^a \in T^*_xM$, 
 the tensors  $\hat J^*\omega|_{\alpha}$  and $g^J_\alpha$ 
have the following components:\par
\medskip
$$\hat J^*\omega|_{\alpha} = J_{i}^a(x) ( d p_a\otimes d x^i - d x^i\otimes d p_a) + 
p_a \left(J^a_{j,i}(x) - J^a_{i,j}(x)\right) d x^i\otimes d x^j \ ,
\eqno(2.7)
$$
\medskip
$$\g^J|_{\alpha} = 
\frac 1 2 p_a N^a_{i\ell}(x) J^\ell_j(x) d x^i\otimes d x^j  = $$
$$ = \frac 1 2 p_a \left[ J^m_i(x) J^a_{\ell, m}(x) -  J^m_\ell(x) J^a_{i, m}(x) - 
J^a_m(x) J^m_{\ell, i}(x) +  \right. $$
$$+ \left. J^a_m(x) J^m_{i, \ell}(x)\right] J^\ell_j(x)
   d x^i\otimes d x^j = $$
   $$ = 
    \frac 1 2 p_a \left\{\left[J^m_i(x) J^\ell_j(x) - J^m_j(x) J^\ell_i(x)\right]J^a_{\ell, m}(x)  +  \left(J^a_{i, j}(x)  - J^a_{j,i}(x)\right)  \right\}   d x^i\otimes d x^j \ .
\eqno(2.8)
$$
\medskip
In the above formulae, we denoted by   $J^i_j(x)$ and $N^a_{i\ell}(x)$ the component of $J_x$ and $N^J_x$  in the coordinate frames of $M$ and by  
  $J^a_{i,j}(x)$
 the 
partial derivatives  $J^a_{i,j}(x) \= \left.\frac{\partial J^a_i}{\partial x^j}\right|_x$. \par
Since
$\omega^{-1}|_{\alpha} =
\frac{\partial}{\partial x^a} \otimes \frac{\partial}{\partial p_a}
- \frac{\partial}{\partial p_a} \otimes \frac{\partial}{\partial x^a}$, we immediately obtain 
that $\J|_{\alpha}$ is of the form
$$\J|_{\alpha} =
J_{i}^a(x) dx^i \otimes \frac{\partial}{\partial x^a} +  J_{i}^a(x)  d p_a \otimes \frac{\partial}{\partial p_i} +  $$
$$ + 
 \frac 1 2 p_a \left\{ \left[J^m_i(x) J^\ell_j(x) - J^m_j(x) J^\ell_i(x)\right]J^a_{\ell, m}(x)  -  \left[J^a_{i, j}(x)  - J^a_{j,i}(x)\right]  \right\}    d x^i\otimes \frac{\partial}{\partial p_i} \ .\eqno(2.9)$$
From this explicit expression, one can directly check that requirement a) holds. Furthermore, 
it is also clear that {\it if $V\in T(T^*M)$ is a vertical vector, then also the vector $\J(V)$ is  vertical\/}. 
\bigskip
\bigskip
\subsection{Hypersurfaces in almost complex manifolds  and their Levi forms} For a submanifold  $S \subset M$ of an almost complex manifold $(M,J)$, we call 
{\it J-invariant} (or {\it $J$-holomorphic\/}) {\it distribution of $S$\/}  
the   family of subspace $\D_x\subset T_xS$, $x\in S$, 
defined by
$$\D_x = \{\ v\in T_x S\ :\ J(v) \in T_xS\ \}\ .\eqno(2.10)$$ 
\medskip
\begin{definition}{\rm A submanifold $S \subset M$ in an almost complex
manifold $(M,J)$ is called {\it totally real\/} if the $J$-invariant 
subspaces (2.10) are trivial at any point}. 
\end{definition}  
\medskip
Assume now that $\Gamma$ is a hypersurface in $M$. 
Notice that, in case $\Gamma$ is (locally) defined as zero set of a smooth real valued function $\rho$
(i.e. 
$\Gamma = \{\ x\in M\ :\ \rho(x) = 0\ \}$),  we may consider the 1-form  
$$\vartheta_x = \left(d\rho\circ J\right)|_{T_x\Gamma}\ .\eqno(2.11)$$
Such 1-form satisfies
$$\ker \vartheta|_x = \D_x\ \qquad \text{for any}\ x\in \Gamma\ .\eqno(2.12)$$
We call   any 1-form $\vartheta$ satisfying (2.12)  a {\it defining 1-form for $\D$\/}. 
The {\it Levi form at $x\in \Gamma$ 
associated with a defining 1-form $\vartheta$\/}  is the quadratic form 
$$\L_x(v) \= - d\vartheta_x(v, Jv)\ ,\qquad \text{for any}\ v\in T_x \Gamma\ .\eqno(2.13)$$
 Notice that, for any vector field $X^{(v)}$ with values in $\D$
such that $X^{(v)}|_x = v$, we may write 
$$\L_x(v) = - d\vartheta_x(X^{(v)}, J X^{(v)}) = $$
$$ = - X^{(v)}(\vartheta(JX^{(v)}))|_x +  JX^{(v)}(\vartheta(X^{(v)}))|_x + \vartheta([X^{(v)}, JX^{(v)}])|_x = $$
$$ = \vartheta([X^{(v)}, JX^{(v)}]) |_x\ ,$$
where we used the fact that, by construction,  $\vartheta(X^{(v)}) \equiv \vartheta(J X^{(v)}) \equiv 0$.\par
The above identity shows that, in case of an integrable complex structures $J$,   the Levi form  $\L_x$ defined in (2.13) coincides with the 
Levi form as classically defined in the theory of functions of several complex variables. \par
\medskip
As in the classical case,  up to multiplication 
by a non-zero real number,  the Levi form
$\L_x$ does not depend on the choice of the defining 1-form $\vartheta$. Moreover, by polarization,  we may say that 
$\L_x$ is the quadratic form associated with the symmetrized bilinear form
$(\BL)^s_x$, where $\BL_x$ is 
$$\BL_x: \D_x\times \D_x \to \R\ ,\qquad \BL_x(v,w) = - d\vartheta_x(v, Jw)$$
and $(\BL)^s_x$ is 
the symmetric part of $\BL_x$, i.e.
$$(\BL_x)^s(v,w) = - \frac{1}{2} \left(d\vartheta_x(v, Jw) + d\vartheta_x(w, Jv)\right) \ .$$
If $J$ is an integrable  complex structure (i.e. $ N^J \equiv 0$), it is not difficult to check   that  the bilinear form $\BL_x$ is symmetric and  $J$-invariant 
and  that  $\L_x$ is the quadratic form associated with  $\BL_x$.  \par
\medskip
\begin{definition}{\rm We  say that a hypersurface $\Gamma\subset M$ is {\it strongly pseudoconvex\/} 
if, for some choice of the defining 1-form $\vartheta$,   the Levi form $\L_x$ is positive definite at any point $x\in \Gamma$.\par
We say that $\Gamma$ is {\it Levi non-degenerate\/}
 if, for some choice of the defining 1-form $\vartheta$, 
the symmetric form $(\BL_x)^s$ is non-degenerate at any point $x\in \Gamma$. }
\end{definition}
\medskip
It is clear that if $\Gamma$ is strongly pseudoconvex (resp. Levi non-degenerate), for {\it any\/} choice of the defining form $\vartheta$, 
the corresponding Levi form $\L_x$ is either positive definite or negative definite (resp. the 
symmetric form $(\BL_x)^s$ is non-degenerate).\par
\medskip
\begin{rem}{\rm In case  $J$ is an integrable complex structure, the Levi non-degeneracy 
condition is equivalent to the non-degeneracy  of $d\vartheta|_{\D\times \D}$, 
i.e. to the hypothesis  that {\it $\D$ is a contact distribution\/}.  On the other hand, 
the reader should be aware that such equivalence 
is no longer valid if $J$ is not integrable, because, in general, 
 the non-degeneracy of $(\BL_x)^s$ does not give
any information on the non-degeneracy of the bilinear form 
$\BL$ (and hence on the non-degeneracy of  $d\vartheta|_{\D\times \D}$). 
If $J$ is not integrable, it is possible to infer that $\D$ is contact, only if  $\Gamma$ is strongly pseudoconvex or under more explicit  conditions on $\BL$.}
\end{rem}
\bigskip
\bigskip
\section{Total reality of the conormal bundles of strongly pseudoconvex hypersurfaces}
\bigskip
In what follows,  $(M,J)$ will denote an almost complex manifold of 
real dimension $2n$, 
$\Gamma$ a hypersurface in $M$, endowed with the $J$-invariant distribution $\D$, and 
$\N^*(\Gamma) \subset T^*M$ the conormal bundle of $\Gamma$ (see definition 
in the Introduction). We will also denote by 
$N$ the Nijenhuis tensor of $J$, by 
$\J$ the lifted almost complex structure of $T^*M$ and by 
$\hat \D$ the 
$\J$-invariant distribution of the submanifold $\N^*(\Gamma)$.\par
\medskip
\begin{lem} \label{lemma1} For any $\alpha \in \N^*(\Gamma)$, the projection $\pi_*: T_\alpha(\N^*(\Gamma)) 
\to T_{\pi(\alpha)} \Gamma$ maps injectively
the $\J$-invariant subspace 
$\hat \D_\alpha \subset T_{\alpha}( \N^*(\Gamma))$ 
onto a  subspace of the J-invariant subspace $\D_{\pi(\alpha)} \subset T_{\pi(\alpha)} \Gamma$.
\end{lem}
\begin{pf} It suffices to show that the vertical 
subspace 
$$\V_\alpha 
 = \{\  V\in T_\alpha(\N^*(\Gamma))\ :\ \pi_*(V) = 0\ \} $$
has trivial intersection with $\hat \D_\alpha$. But this is a direct consequence
of the following two facts that: 1)   $\J$ maps vertical vectors 
into vertical vectors and hence  $\V_\alpha \cap \hat \D_\alpha$ is  an even dimensional
$\J$-invariant subspace of $\V_\alpha$; 2) $\V_\alpha$ has 
 dimension one. 
\end{pf}
From the previous lemma,  we have that for any $0 \neq V\in \hat \D_\alpha$, the projection 
$\pi_*(V)$ must be  a non-trivial element of $\D_{\pi(\alpha)}$. \par
\medskip
\begin{lem} \label{lemma2} $\N^*(\Gamma)$ is a  Lagrangian submanifold of 
$(T^*M, \omega)$, i.e.  for any $\alpha \in \N^*(\Gamma)$ and any
$V, W \in T_\alpha(\N^*(\Gamma))$
$$\omega_\alpha(V,W) = 0\ .$$
\end{lem}
\begin{pf} Consider two vector fields $X^{(V)}, X^{(W)}$ in $\N^*(\Gamma)$ such that $X^{(V)}|_\alpha = V$ 
and $X^{(W)}|_\alpha = W$. By definitions, 
$$\omega_\alpha(V,W) = d\theta_\alpha(X^{(V)},X^{(W)}) = $$
$$ = X^{(V)}(\theta(X^{(W)}))|_\alpha - X^{(W)}(\theta(X^{(V)}))|_\alpha 
- \theta_\alpha([X^{(V)},X^{(W)}]) = 0\ ,$$
because  for any vector field $Z$ on $\N^*(\Gamma)$, 
$\pi_*(Z) \in T\Gamma$ and hence 
$\theta_\beta(Z) = \beta(\pi_*(Z)) = 0$ for any
$\beta \in \N^*(\Gamma)$. 
\end{pf}
\bigskip
Now, we can state and prove our main theorem.\par
\medskip
\begin{theo} \label{maintheorem}ÊLet $\Gamma \subset M$ be a hypersurface in 
an almost complex manifold $(M,J)$ and assume that at least one of the following 
hypothesis  is satisfied:
\begin{itemize}
\item[i)] $J$ is integrable and the $J$-invariant distribution is contact (i.e. $\Gamma$ is Levi non-degenerate); 
\item[ii)] $\Gamma$ is strongly pseudo-convex. 
\end{itemize} 
Then  $\N^*(\Gamma) \setminus\{ \operatorname{zero\ section}\}$ 
is a totally real submanifold of $T^*M$. 
\end{theo}
\begin{pf} Suppose not and assume that, for some $\alpha
 \in \N^*(\Gamma) \setminus\{ \operatorname{zero\ section}\}$, there 
exists a non-trivial vector $V$, which belongs to  the $\J$-invariant 
subspace $\hat \D_\alpha$, i.e.  $0 \neq \J(V) \in T_\alpha \N^*(\Gamma)$.  
By Lemma \ref{lemma2}, we also have that for any $W\in T_\alpha \N^*(\Gamma)$
$$\omega_\alpha(\J(V),W) = 0\ .\eqno(3.1)$$
On the other hand, by (2.2) and the definition of $\J$
$$\omega_\alpha(\J(V), W) = 
d(\hat J^* \theta)_\alpha(V,W) + g^J_\alpha(V,W)\ .\eqno(3.2)$$
Formulae (3.1) and (3.2) imply  that, for any vector fields $X^{(V)}$ and $X^{(W)}$ on $\N^*(\Gamma)$ such that 
$$X^{(V)}|_\alpha = V\ ,\qquad X^{(W)}|_\alpha = W\ ,$$
we have that 
$$X^{(V)}(\hat J^*\theta(X^{(W)}))|_\alpha - 
X^{(W)}(\hat J^*\theta(X^{(V)}))|_\alpha - \hat J^*\theta_\alpha([X^{(V)},X^{(W)}]) = $$
$$ = - \frac 1 2 
\alpha(N(\pi_*(X^{(V)}), J\pi*(X^{(W)}))_{\pi(\alpha) })
\ .\eqno(3.3)$$
Now, recall that,  
by Lemma \ref{lemma1}, $ v \= \pi_*(V) \neq 0$ and belongs to $\D_{\pi(\alpha)}$. 
Assume now that  the vector $W$ is 
so  that $w \= \pi_*(W) \in \D_{\pi(\alpha)}$ and 
 let us denote by $\tilde X^{(V)}$ and $\tilde X^{(W)}$ the two vector fields in $T\Gamma$
defined by  
$$\tilde X^{(V)} \= \pi_*(X^{(V)})\ ,\qquad \tilde X^{(W)}\= \pi_*(X^{(W)})\ .$$  
 Clearly, $ \tilde X^{(V)}|_{\pi(\alpha)} = v$ and $ \tilde X^{(W)}|_{\pi(\alpha)} = w$.\par
 Now, notice that if $\rho$ is a local defining function for $\Gamma$, then any $\beta 
 \in \N^*(\Gamma) \setminus \{\text{zero section}\}$ is of the form
 $\beta = \lambda \cdot d\rho$, for some $\lambda \in \R\setminus\{0\}$.  Hence, 
 by definition of the tautological form $\theta$ and since $\hat J$ is 
fiber preserving, for any 
$\beta \in \N^*(\Gamma)\setminus\{ \operatorname{zero\ section}\}$ 
and any vector field on $\N^*(\Gamma)$
$$\hat J^*\theta_\beta(Z) = \beta(J \pi_*(Z)) = \lambda \cdot d\rho(J \pi_*(Z)) = 
\lambda \cdot \vartheta(\pi_*(Z))$$
where $\vartheta = d\rho \circ J$  and $\lambda \in \R\setminus \{0\}$.  Recall that $\vartheta$ satisfies (2.11). Now, if 
$\alpha = \lambda_o \cdot d\rho|_{\pi(\alpha)}$ for some fixed $\lambda_o \neq 0$, then 
the left hand side of 
(3.3) can be written as 
$$\left.X^{(V)}\left(\lambda \cdot \vartheta(\tilde X^{(W)})\right)\right|_\alpha - 
X^{(W)}\left.\left(\lambda \cdot \vartheta(\tilde X^{(V)})\right)\right|_\alpha 
- \lambda_o \cdot \vartheta_{\pi(\alpha)}([\tilde X^{(V)}, \tilde X^{(W)}]) = $$
$$ = \left.X^{(V)}(\lambda) \right|_\alpha \cdot \vartheta_{\pi(\alpha)}(w) - 
\left.X^{(W)}(\lambda) \right|_\alpha \cdot
\vartheta_{\pi(\alpha)}( v) + $$
$$ + \lambda_o\cdot \left(\tilde  X^{(V)}(\vartheta(\tilde  X^{(W)}))|_{\pi(\alpha)} - 
\tilde X^{(W)}(\vartheta(\tilde  X^{(V)}))|_{\pi(\alpha)} - 
\vartheta_{\pi(\alpha)}([\tilde X^{(V)}, \tilde  X^{(W)}])\right) = $$
$$ =  \left.X^{(V)}(\lambda)\right|_\alpha\cdot \vartheta_{\pi(\alpha)}(w) - 
\left. X^{(W)}(\lambda)\right|_\alpha  \cdot
\vartheta_{\pi(\alpha)}(v) + 
 \lambda_o\cdot  d\vartheta_{\pi(\alpha)}(v, w)\ .
\eqno(3.4)
$$
Since we are assuming that $v$ and $ w$ are both in $\D$, we have that 
$\vartheta_{\pi(\alpha)}( v) = \vartheta_{\pi(\alpha)}(w) = 0$ and hence (3.3) 
reduces to
$$\lambda_o \cdot d\vartheta_{\pi(\alpha)}(v, w) =  - \frac 1 2 
\alpha(N(v, J w))\ .\eqno(3.5)$$
This should be true for any choice of 
$w \in \D_{\pi(\alpha)}$.  This 
gives a contradiction if $N \equiv 0$ (i.e. $J$ is integrable) 
and $\D$ is contact, because in this case there should exist
a vector $w\in \D_{\pi(\alpha)}$ so that  $d\vartheta_{\pi(\alpha)}(v, w) \neq 0$. 
In case  $\Gamma$ is strongly pseudo-convex and
$N$ is not necessarily  equal to $0$,  in order to  
get a contradiction  it suffices  to assume that $ w =  J v$. In fact,
in this case, $N(v, J w) = - N(v, v) = 0$, while, 
by strong pseudoconvexity $
\lambda_o \cdot d\vartheta_{\pi(\alpha)}(v,  J v) = - \lambda_o
\cdot \L_{\pi(\alpha)}( v) \neq 0$. 
\end{pf}

\bigskip
\bigskip
\font\smallsmc = cmcsc8
\font\smalltt = cmtt8
\font\smallit = cmti8
\hbox{\parindent=0pt\parskip=0pt
\vbox{\baselineskip 9.5 pt \hsize=3.1truein
\obeylines
{\smallsmc
Andrea Spiro
Dip. Matematica e Informatica
Universit\`a di Camerino
Via Madonna delle Carceri
I-62032 Camerino (Macerata)
ITALY
}\medskip
{\smallit E-mail}\/: {\smalltt andrea.spiro@unicam.it}
}
}

\end{document}